\documentclass[12pt,reqno]{amsart}
\input{avv.sty}
\title{Generalized Convexity and Inequalities}
\author{G. D. Anderson, M. K. Vamanamurthy, and M. Vuorinen}
\date{}
\subjclass[2000]{Primary  33C05, 33C20. Secondary  26A51.}

\keywords{ Convexity, monotonicity, power series, hypergeometric function, generalized hypergeometric series.}

\usepackage{amssymb}
\usepackage{amsmath}
\usepackage{amsthm}
\usepackage{setspace} 
\usepackage{geometry}
\newtheorem{lemma}{Lemma}[section]
\newtheorem{theorem}[lemma]{Theorem}
\newtheorem{corollary}[lemma]{Corollary}
\newtheorem{remark}[lemma]{Remark}
\newtheorem{ex}[lemma]{Examples}
\newtheorem{definition}[lemma]{Definition}
\newtheorem{notation}[lemma]{Notation}
\newtheorem{Pf13}[lemma]{Proof of Theorem 1.3}
\newtheorem{Pf14}[lemma]{Proof of Theorem 1.4}
\newtheorem{Pf15}[lemma]{Proof of Theorem 1.5}
\renewcommand{\theequation}{\thesection.\arabic{equation}}

\begin{document}
\markboth{\textsc{G. D. ANDERSON, M. K. VAMANAMURTHY, AND
M. VUORINEN}}{\textsc{}}

\begin{abstract} 
  Let $\mathcal{R}_+ =(0,\infty)$ and let $\mathcal{M}$ be the family of all
   mean values of two numbers in $\mathcal{R}_+$ (some examples are the arithmetic,
   geometric, and harmonic means). Given $m_1, m_2 \in {\mathcal{M}},$
   we say that a function $f: \mathcal{R}_+ \to \mathcal{R}_+$ is $(m_1,m_2) $-convex if
   $f(m_1(x,y) ) \le m_2(f(x), f(y))$ for all $x,y \in \mathcal{R}_+ .$
   The usual convexity is the special case when both mean values
are arithmetic means. We study the dependence of 
$(m_1,m_2)$-convexity on $m_1$ and $m_2$ and
give sufficient conditions  for $(m_1, m_2)$-convexity of functions defined by Maclaurin series. The criteria involve
the Maclaurin coefficients. Our results yield a class of new
   inequalities for several special functions such as the Gaussian
hypergeometric function and a generalized Bessel function.
\end{abstract}

\maketitle


\renewcommand{\theequation}{\thesection.\arabic{equation}}

\font\fFt=eusm10 
\font\fFa=eusm7  
\font\fFp=eusm5  
\def\K{\mathchoice
{\hbox{\,\fFt K}}
{\hbox{\,\fFt K}}
{\hbox{\,\fFa K}}
{\hbox{\,\fFp K}}}

\def\E{\mathchoice
{\hbox{\,\fFt E}}
{\hbox{\,\fFt E}}
{\hbox{\,\fFa E}}
{\hbox{\,\fFp e}}}

\font\fFt=eusm10 
\font\fFa=eusm7  
\font\fFp=eusm5  
\def\A{\mathchoice
{\hbox{\,\fFt A}}
{\hbox{\,\fFt A}}
{\hbox{\,\fFa A}}
{\hbox{\,\fFp A}}}

\section{Introduction}

In this paper we study several convexity and monotonicity properties
of certain functions and deduce sharp inequalities. We deduce
analogous results for certain power series, especially
hypergeometric functions.  This work continues studies in
\cite{2}, \cite{7}, and \cite{8}.

        The following result \cite[Theorem 4.3]{12}, a variant of a result by Biernacki and Krzy\.z \cite{9}, will be very useful in studying convexity and monotonicity of certain power series.
\vspace{.2in}

\begin{lemma}\label{lem1}
For $0 < R \le \infty$, let $f(x) = \sum_{n=0}^{\infty} a_nx^n$ and
$g(x) = \sum_{n=0}^{\infty} b_nx^n$ be two real power series converging on
the interval $(-R,R)$.  If the sequence $\{a_n/b_n\}$ is increasing
(decreasing), and $b_n > 0$ for all $n$, then the function
$$
h(x) = \frac{f(x)}{g(x)} = \dfrac{\sum_{n=0}^{\infty} a_nx^n}{\sum_{n=0}^{\infty} b_nx^n}
$$
is also increasing (decreasing) on $(0,R)$.  In fact, the function
$$
g(x)f'(x)-f(x)g'(x)
$$ 
has positive Maclaurin coefficients.
\end{lemma}
\vspace{.2in}

\begin{notation}  If  $f(x) = \sum_{n=0}^{\infty }a_nx^n$ and $g(x) =\sum_{n=0}^{\infty }b_nx^n$ are two power series, where $b_n > 0$ for all $n$, we let $T_n = T_n(f(x),g(x)) = a_n/b_n$.  We will use $F = F(a,b;c;x)$ to denote the Gaussian hypergeometric function
$$
F(a,b;c;x) = {}_2F_1(a,b;c;x)=\sum_{n=0}^{\infty}\dfrac{(a,n)(b,n)}{(c,n)n!}x^n,\ \ |x|<1,
$$
where  $(a,n)$ denotes the product $a(a+1)(a+2)\cdots (a+n-1)$ when $n\ge 1$, and $(a,0) = 1$ if  $a\ne 0$.  The expression $(0,0)$ is not defined.  
\end{notation}
\vspace{.2in}

Throughout this paper, for $x\in (0,1)$ we denote $x'=\sqrt{1-x^2}$.
\vspace{.2in}

Our main results, to be proved in Section 3, are as follows.
\vspace{.2in}

\begin{theorem}\label{th3}
Let $F(x) = F(a,b;c;x)$, for $a,b,c > 0$ and $|x| < 1$. Then the following results hold.
\begin{enumerate}

\item[(1)] If $ab/(a+b+1) < c$, then $\log F(x)$ is convex on $(0,1)$. In particular,
$F((x+y)/2) \le \sqrt{F(x)F(y)}$, for all $x, y \in (0,1)$, with equality if and only if $x = y$.
\item[(2)] If $(a-c)(b-c) > 0$, then $\log F(1-e^{-t})$ is concave on $(0,\infty)$. In particular, $\sqrt{F(x)F(y)} \le F(1-\sqrt{(1-x)(1-y)} )$, for all $x, y \in (0,1)$, with equality if and only if $x = y$.
\item[(3)] If $a+b \ge c$, then $F(1-e^{-t})$ is convex on $(0,\infty )$. In particular,
$F(1-\sqrt{(1-x)(1-y)}) \le (F(x)+F(y))/2$ for all $x, y \in (0,1)$, with equality iff  $x=y$.
\end{enumerate}
\end{theorem}
\vspace{.2in}

\begin{theorem}\label{newth2}
Let $F(x) = F(a,b;c;x)$, with $a,b,c > 0$ and $|x| < 1$. If
$a+b \ge c \ge 2ab$ and $c > a+b -1/2$, then $1/F(x)$ is concave on $(0,\infty )$. In particular,
$$
F\left(\frac{x+y}2\right) \le \frac{2F(x)F(y)}{F(x)+F(y)},
$$
for all $x, y \in (0,1)$, with equality if and only if $x=y$.
\end{theorem}
\vspace{.2in}

\begin{theorem}\label{mf}
For $0 < R < \infty$, let $f(x) = \sum_{n= 0}^{\infty }a_n x^n$, $a_n > 0$, be convergent on $(-R, R)$.
Let $m = m_f$ be the function  defined by $m(x) = f(R-x^2/R)/f(x^2/R)$. If the sequence $\{R(n+1)a_{n+1}/a_n -n\}$  is decreasing, then
$$
\frac{1}{m(\sqrt[4]{(R^2-x^2)(R^2-y^2)})} \le \sqrt{m(x)m(y)} \le m(\sqrt{xy})
$$
for all $x,y \in (0,R)$, with equality if and only if $x=y$.
\end{theorem}
\vspace{.2in}

The hypergeometric function contains, as its limiting or
special cases, many well-known special functions.
For instance, the reader may find in \cite{15} a list
that contains several hundreds of rational triples
$(a,b,c)$ such that the hypergeometric function
$F(a,b;c;r)$ reduces to a well-known function.
Therefore, the above results yield new inequalities even for many familiar elementary functions.
\vspace{.2in}

\section{Generalized convexity}

The notions of convexity and concavity of a real-valued function of a real
variable are well known \cite{17}. In this section we study certain
generalizations of these notions for a positive-valued function of a
positive variable.
\vspace{.2in}

\begin{definition}
A function $M : (0,\infty )\times (0,\infty ) \to (0,\infty )$ is
called a {\it Mean function} if 
\begin{enumerate}

\item[(1)] $M(x,y) = M(y,x)$, 
\item[(2)] $M(x,x) =
x$,
\item[(3)] $x < M(x,y) < y$, whenever $x < y$,
\item[(4)] $M(ax,ay)=aM(x,y)$ for all $a>0$.
\end{enumerate}
\end{definition}
\vspace{.2in}

\begin{ex}\cite{10},\cite{11}
\begin{enumerate}

\item[(1)] $M(x,y) = A(x,y) = (x+y)/2$ is the {\it Arithmetic Mean}.
\item[(2)] $M(x,y) = G(x,y) = \sqrt {xy}$ is the {\it Geometric Mean}.
\item[(3)] $M(x,y) = H(x,y) = 1/A(1/x,1/y)$ is the {\it Harmonic Mean}.
\item[(4)] $M(x,y) = L(x,y) = (x-y)/(\log x - \log y)$ for $x\ne y$, and $L(x,x) = x$, is the {\it Logarithmic Mean}.
\item[(5)] $M(x,y) = I(x,y) = (1/e)(x^x/y^y)^{1/(x-y)}$ for $x\ne y$, and $I(x,x) = x$, is the {\it Identric Mean}.
\end{enumerate}
\end{ex}
\vspace{.2in}

\begin{definition}
Let $f : I \to (0,\infty )$ be continuous, where $I$ is a
subinterval of $(0,\infty )$. Let $M$ and $N$ be any two Mean
functions. We say $f$ is $MN$-convex (concave) if $f(M(x,y)) \le
(\ge )\ N(f(x),f(y))$, for all $x,y\in I$.
\end{definition}
\vspace{.2in}

Note that this definition reduces to usual convexity (concavity)
when $M = N = A$. The concept of $MN$-convexity has been studied 
extensively in the literature from various points of view
(see e.g. \cite{5}, \cite{1}, \cite{13}, \cite{14}), but so far as we know, very few criteria for $MN$-convexity of Maclaurin series are known in terms of the coefficients.  We will concentrate on criteria of this type in the main results of this paper in Section 3. We now show that for $M, N = A, G, H$, the nine
possible $MN$-convexity (concavity) properties reduce to ordinary
convexity (concavity) by a simple change of variable.
\vspace{.2in}

\begin{theorem}\label{gen}
Let $I$ be an open subinterval of $(0,\infty)$ and let $f \colon I
\to (0,\infty)$ be continuous. In parts $(4)-(9)$, let $I =(0,b)$, $0 < b < \infty$.
\begin{enumerate}

\item[(1)] $f$ is $AA$-convex (concave) if and only if  $f$ is convex (concave).
\item[(2)] $f$ is $AG$-convex (concave) if and only if  $\log f$ is convex (concave).
\item[(3)] $f$ is $AH$-convex (concave) if and only if $1/f$ is concave (convex).
\item[(4)] $f$ is $GA$-convex (concave) on $I$ if and only if $f(b e^{-t})$ is convex (concave) on $(0,\infty )$.
\item[(5)] $f$ is $GG$-convex (concave) on $I$ if and only if $\log f(b e^{-t})$ is convex (concave) on $(0,\infty )$.
\item[(6)] $f$ is $GH$-convex (concave) on $I$ if and only if $1/f(b e^{-t})$ is concave (convex) on $(0,\infty )$.
\item[(7)] $f$ is $HA$-convex (concave) on $I$ if and only if $f(1/x)$ is convex (concave) on $(1/b,\infty )$.
\item[(8)] $f$ is $HG$-convex (concave) on $I$ if and only if $\log f(1/x)$ is convex (concave) on $(1/b,\infty )$.
\item[(9)] $f$ is $HH$-convex (concave) on $I$ if and only if $1/f(1/x)$ is concave (convex) on $(1/b,\infty )$.
\end{enumerate}
\end{theorem}

\begin{proof}
\begin{enumerate}

\item[(1)] This follows by definition.
\item[(2)]
\begin{align*}
f(A(x,y)) &\le (\ge )\ G(f(x),f(y))\\
 &\iff f((x+y)/2) \le (\ge )\ \sqrt {f(x)f(y)}\\
 &\iff \log f((x+y)/2) \le (\ge )\ \frac{1}{2}\left(\log f(x) + \log f(y)\right),
\end{align*}
\noindent hence the result.
\item[(3)]
\begin{align*}
f(A(x,y)) &\le (\ge )\ H(f(x),f(y))\\
 &\iff f((x+y)/2) \le (\ge )\ 2/(1/f(x) + 1/f(y))\\
&\iff 1/f\left((x+y)/2\right) \ge (\le )\ \frac{1}{2}(1/f(x) + 1/f(y)),
\end{align*}
hence the result.
\item[(4)] With $x = b e^{-r}$ and $y = b e^{-s}$,
\begin{align*}
f(G(x,y))  &\le (\ge )\ A(f(x)+f(y))\\
&\iff f\left(b e^{-(r+s)/2}\right) \le (\ge )\ \frac{1}{2}\left(f(b e^{-r}) + f(b e^{-s})\right),
\end{align*}
\noindent hence the result.
\item[(5)] With $x = b e^{-r}$ and $y = b e^{-s}$,
\begin{align*}
&f(G(x,y)) \le (\ge )\ G(f(x),f(y))\\
&\iff \log f(b e^{-(r+s)/2}) \le (\ge )\ \frac{1}{2}(\log f(b e^{-r}) + \log f(b e^{-s})),
\end{align*}
\noindent hence the result.
\item[(6)] With $x = b e^{-r}$ and $y = b e^{-s}$,
\begin{align*}
f(G(x,y)) &\le (\ge )\ H(f(x),f(y))\\
 &\iff 1/f(b e^{-(r+s)/2}) \ge (\le )\ \frac{1}{2}(1/f(b e^{-r}) + 1/f(b e^{-s})),
\end{align*}
hence the result.
\item[(7)] Let $g(x) = f(1/x)$, and let $x, y \in (1/b,\infty )$, so that $1/x, 1/y \in (0,b)$. Then
$f$ is $HA$-convex (concave) on $(0,b)$ if and only if
\begin{align*}
f\left(\frac{2}{x+y}\right) &\le (\ge )\ (1/2)(f(1/x)+f(1/y))\\
&\iff g\left(\frac{x+y}{2}\right) \le (\ge )\ \frac{1}{2}(g(x)+g(y)),
\end{align*}
hence the result.
\item[(8)] Let $g(x) = \log f(1/x)$, and let $x, y \in (1/b,\infty )$, so $1/x,1/y \in (0,b)$. Then $f$ is $HG$-convex (concave) on $(0,b)$ if and only if
\begin{align*}
&f\left(\frac{2}{x+y}\right) \le (\ge )\ \sqrt {f(1/x)f(1/y)}\\
&\iff \log f((2/(x+y))) \le (\ge )\ (1/2)(\log f(1/x)+\log f(1/y))\\
 &\iff g\left(\frac{x+y}{2}\right) \le (\ge )\ \frac{1}{2}(g(x)+g(y)),
\end{align*}
hence the result.
\item[(9)] Let $g(x) = 1/f(1/x)$, and let $x, y \in (1/b,\infty )$, so $1/x,1/y \in (0,b)$. Then $f$ is $HH$-convex (concave) on $(0,b)$ if and only if
\begin{align*}
&f\left(\frac{2}{x+y}\right) \le (\ge )\ 2/(1/f(1/x)+1/f(1/y))\\
&\iff 1/f(2/(x+y)) \ge (\le )\ \frac{1}{2}(1/f(1/x)+1/f(1/y))\\
&\iff g\left(\frac{x+y}{2}\right) \ge (\le )\ \frac{1}{2}(g(x)+g(y)),
\end{align*}
hence the result.
\end{enumerate}
\end{proof}
\vspace{.2in}

The next result is an immediate consequence of Theorem \ref{gen}.
\vspace{.2in}

\begin{corollary}\label{gencor}
Let $I$ be an open subinterval of $(0,\infty)$ and let $f \colon I
\to (0,\infty)$ be differentiable. In parts $(4)-(9)$, let $I =
(0,b)$, $0 < b < \infty$.
\begin{enumerate}

\item[(1)] $f$ is $AA$-convex (concave) if and only if $f'(x)$ is increasing (decreasing).
\item[(2)] $f$ is $AG$-convex (concave) if and only if $f'(x)/f(x)$ is increasing (decreasing).
\item[(3)] $f$ is $AH$-convex (concave) if and only if $f'(x)/f(x)^2$ is increasing (decreasing).
\item[(4)] $f$ is $GA$-convex (concave) if and only if $xf'(x)$ is increasing (decreasing).
\item[(5)] $f$ is $GG$-convex (concave) if and only if $xf'(x)/f(x)$ is increasing (decreasing).
\item[(6)] $f$ is $GH$-convex (concave) if and only if $xf'(x)/f(x)^2$ is increasing (decreasing).
\item[(7)] $f$ is $HA$-convex (concave) if and only if $x^2 f'(x)$ is increasing (decreasing).
\item[(8)] $f$ is $HG$-convex (concave) if and only if $x^2 f'(x)/f(x)$ is increasing (decreasing).
\item[(9)] $f$ is $HH$-convex (concave) if and only if $x^2 f'(x)/f(x)^2$ is increasing (decreasing).
\end{enumerate}
\end{corollary}
\vspace{.2in}

\begin{remark}
Since $H(x,y) \le G(x,y) \le A(x,y)$, it follows that
\begin{enumerate}

\item[(1)] $f$ is $AH$-convex $\Longrightarrow$ $f$ is $AG$-convex $\Longrightarrow$ $f$ is $AA$-convex.
\item[(2)] $f$ is $GH$-convex $\Longrightarrow$ $f$ is $GG$-convex $\Longrightarrow$ $f$ is $GA$-convex.
\item[(3)] $f$ is $HH$-convex $\Longrightarrow$ $f$ is $HG$-convex $\Longrightarrow$ $f$ is $HA$-convex.
\end{enumerate}
Further, if f is increasing (decreasing) then $AN$-convex (concave)
implies $GN$-convex (concave) implies $HN$-convex (concave), where
$N$ is any Mean function. For concavity, the implications in (1), (2), and (3) are reversed.
These implications are strict, as shown by the examples below.
\end{remark}
\vspace{.2in}

\begin{ex}
\begin{enumerate}

\item[(1)] $f(x) = \cosh x$ is $AG$-convex, hence $GG$-convex and $HG$-convex, on $(0,\infty )$. But it is not $AH$-convex, nor $GH$-convex, nor $HH$-convex.
\item[(2)] $f(x) = \sinh x$ is $AA$-convex, but $AG$-concave on $(0,\infty )$.
\item[(3)] $f(x) = e^x$ is $GG$-convex and $HG$-convex, but neither $GH$-convex nor $HH$-convex, on $(0,\infty )$.
\item[(4)] $f(x) = \log (1+x)$ is $GA$-convex, but $GG$-concave on $(0,\infty )$.
\item[(5)] $f(x)= \arctan x$ is $HA$-convex, but not $HG$-convex, on $(0,\infty )$.
\end{enumerate}
\end{ex}
\vspace{.2in}

\section{Applications to power series}

\begin{theorem}\label{th2}
Let $f(x) = \sum_{n=0}^{\infty }a_n x^n$, where $a_n > 0$ for $n= 0, 1, 2,\ldots$, be convergent on $(-R,R)$, $0<R<\infty $. Then the following convexity results hold.
\begin{enumerate}

\item[(1)] $f$ is $AA$-convex and $GG$-convex.
\item[(2)] If the sequence $\{(n+1)a_{n+1}/a_n\}$ is increasing (decreasing), then $f$ is $AG$-convex (concave) on $(0,R)$.  In particular,
$$
f\left(\frac{x+y}2\right) \le (\ge )\ \sqrt{f(x)f(y)}
$$
for all $x, y \in (0,R)$, with equality if and only if  $x = y$.
\item[(3)] Let $b_n = \sum_{k = 0}^n a_k a_{n-k}$. If the sequence $\{(n+1)a_{n+1}/b_n\}$ is increasing (decreasing), then $f$ is $AH$-convex (concave) on $(0,R)$.
\item[(4)] Let $b_n = \sum_{k = 0}^n a_k a_{n-k}$. If the sequence $\{na_n/b_n\}$ is increasing (decreasing), then $f$ is $GH$-convex (concave) on $(0,R)$.
\item[(5)] If the sequence $\{R(n+1)a_{n+1}/a_n - n\}$  is increasing (decreasing), then the function $(R-x)f'(x)/f(x)$ is increasing (decreasing) on $(0,R)$, so that the function $\log f(R(1-e^{-t}))$ is convex (concave) on $(0,\infty )$. In particular,
$$
f(R-\sqrt{(R-x)(R-y)}\ )\le (\ge )\ \sqrt{f(x)f(y)}  
$$
for all $x, y \in (0,R)$, with equality if and only if $x = y$.
\item[(6)] If the sequence $\{na_nR^n\}$ is increasing (decreasing), then the function $(R-x)f'(x)$ is increasing (decreasing) on $(0,R)$, so that the function $f(R(1-e^{-t}))$ is convex (concave) as a function of $t$ on $(0,\infty )$. In particular,
$$
f(R-\sqrt{(R-x)(R-y)}\ ) \le (\ge )\ \frac{f(x)+f(y)}{2},
$$
for all $x, y \in (0,R)$, with equality if and only if $x = y$.
\item[(7)] If the sequence $\{na_nR^n\}$ is increasing and if also the sequence $\{n!a_nR^n/(1/2,n)\}$ is decreasing, then the function $1/f(x)$ is concave on $(0,R)$. In particular,
$$
f\left(\frac{x+y}2\right) \le \frac{2f(x)f(y)}{f(x)+f(y)},
$$
for all $x, y \in (0,R)$, with equality if and only if $x=y$.
\end{enumerate}
\end{theorem}
\begin{proof}
\begin{enumerate}

\item[(1)]  These follow trivially from Corollary \ref{gencor} and Lemma \ref{lem1}.

\item[(2)] $T_n(f'(x),f(x))=(n+1)a_{n+1}/a_n$,
which is increasing (decreasing).  Thus, by Lemma \ref{lem1} and Corollary \ref{gencor}(2) the assertion follows.
\item[(3)] Since $T_n(f'(x),f(x)^2) = (n+1)a_{n+1}b_n$, the result follows by Lemma \ref{lem1} and Corollary \ref{gencor}(3).
\item[(4)] Since $T_n(xf'(x),f(x)^2) = na_n/b_n$,
the result follows by Lemma \ref{lem1} and Corollary \ref{gencor}(6),(9).
\item[(5)]
$$
\frac {d}{dt}\log f(R(1-e^{-t})) = Re^{-t} \frac{f'(R(1-e^{-t}))}{f(R(1-e^{-t}))} = (R-x)\frac{f'(x)}{f(x)},
$$
where $x = R(1- e^{-t})$. Then
$$
T_n((R-x)f'(x),f(x)) = R(n+1)a_{n+1}/a_n - n,
$$
which is increasing (decreasing), so that the assertion follows by Lemma \ref{lem1}.
\item[(6)]$(d/dt)f(R(1-e^{-t})) = Re^{-t} f'(R(1-e^{-t})) = (R-x)f'(x) = f'(x)/(1/(R-x))$, where $x = R(1-e^{-t})$. Then, $$
T_n(f'(x),1/(R-x))=(n+1)a_{n+1}R^{n+1},
$$
which is increasing (decreasing) by hypothesis, so that the assertion follows by Lemma \ref{lem1}.
\item[(7)] First,
$$
\frac{d}{dx}\frac{1}{f(x)} =\frac{-f'(x)}{f(x)^2}.
$$
Now $(R-x)f'(x) = f'(x)/[1/(R-x)]$, so that
$$
T_n(f'(x),1/(R-x))=(n+1)a_{n+1}R^{n+1},
$$
which is increasing by hypothesis.  Hence, $(R-x)f'(x)$ is increasing on $(0,R)$, by Lemma \ref{lem1}.  Since
$$
\sqrt{R-x}f(x)=f(x)/(R-x)^{-1/2},
$$
we have 
$$
T_n(f(x),(R-x)^{-1/2})=\frac{n!a_nR^{n+1/2}}{(1/2,n)},
$$
which is decreasing by hypothesis.  Hence, $\sqrt{R-x}f(x)$ is also decreasing on $(0,R)$, by Lemma 1.1. Dividing $(R-x)f'(x)$ by the square of $\sqrt{R-x}f(x)$, we see that  $(d/dx)[1/f(x)]$ is decreasing in $x$ on $(0,R)$, proving the assertion.
\end{enumerate}
\end{proof}
\vspace{.2in}

\begin{Pf13}
\end{Pf13}
\begin{enumerate}

\item[(1)] $T_n(F'(x),F(x)) =  (n+1)a_{n+1}/a_n = (a+n)(b+n)/(c+n)$. Hence,
\begin{align*}
T_{n+1} - T_n > 0 \iff &(a+n+1)(b+n+1)(c+n)\\
&- (a+n)(b+n)(c+n+1) > 0\\
 \iff & (a+n)(c+n) + (b+n)(c+n) +(c+n)\\
 &- (a+n)(b+n) > 0\\
 \iff & n^2 + n(2c+1) + (ac+bc+c-ab) > 0\\
 \iff & ab/(a+b+1) < c.
 \end{align*}
 Hence, the assertion follows from Theorem \ref{th2}(2) and Lemma \ref{lem1}.
\item[(2)]
\begin{align*}
T_n((1-x)F'(x),F(x)) &= \frac{(n+1)a_{n+1}}{a_n} - n\\
&= \frac{(a+n)(b+n)}{c+n} - n\\
&= a+b-c + \frac{(a-c)(b-c)}{c+n},
\end{align*}
which is decreasing if and only if $(a-c)(b-c) > 0$, so that the assertion follows from Theorem \ref{th2}(5).
\item[(3)]
$T_{n-1}(F'(x),1/(1-x)) = n a_n = (a,n)(b,n)/[(c,n)(n-1)!]$. Hence,
\begin{align*}
T_n& - T_{n-1} = \frac{(a,n+1)(b,n+1)}{(c,n+1)n!} - \frac{(a,n)(b,n)}{(c,n)(n-1)!}\\
&= \frac{(a,n)(b,n)}{(c,n+1)n!} [(a+n)(b+n) - n(c+n)]>0\\
&\iff n(a+b-c) + ab > 0\ {\rm for\ all}\ n\\
&\iff a+b \ge c.
\end{align*}
Hence, the assertion follows from Theorem \ref{th2}(6).\hspace{.75in} $\square$
\end{enumerate}
\vspace{.2in}

\begin{theorem}
(cf. \cite[Lemma 2.1]{6}). Let $F(x)  = F(a,b;c;x)$, with $c=a+b, a, b > 0$, and $|x| < 1$. Then
\begin{enumerate}

\item[(1)] $\log F(x)$ is convex on $(0,1)$,
\item[(2)] $\log F(1-e^{-t})$ is concave on $(0,\infty )$,
\item[(3)] $F(1-e^{-t})$ is convex on $(0,\infty )$.
\end{enumerate}
In particular,
$$
F\left(\frac{x+y}2\right) \le \sqrt{F(x)F(y)} \le F(1-\sqrt{(1-x)(1-y)}) \le \frac{F(x)+F(y)}2
$$
for all $x, y \in (0,1)$, with equality if and only if $x = y$.
\end{theorem}
\begin{proof}
This result follows immediately from Theorem \ref{th3}.
\end{proof}
\vspace{.2in}

\begin{theorem}\label{th4}
Let $a,b,c > 0$, $a, b \in (0,1)$ and $a < c, b < c$. Let $F$ and $F_1$ be the conjugate hypergeometric functions on $(0,1)$ defined by
$F = F(x) = F(a,b;c;x)$ and $F_1 = F_1(x) = F(1-x)$. Then the function $f$ defined by
$f(x) = x(1-x)F(x)F_1(x)$ is increasing on $(0,1/2]$ and decreasing on $[1/2,1)$.
\end{theorem}
\begin{proof}
Since $f(x) = f(1-x)$, it is enough to prove the assertion on $(0,1/2]$.
Following Rainville \cite[p. 51]{16} we let $F(a-)=F(a-1,b;c;x)$ and $F_1(a-) = F_1(a-1,b;c;x)$.  Now, since
$$
x(1-x)F'(x) = (c-a)F(a-)+(a-c+bx)F
$$
\cite[Theorem 3.12(2)]{3}, we have
\begin{align*}
f'(x) &= x(1-x)[F'(x)F_1(x) + F(x) F'_1(x)] +(1-2x)F(x)F_1(x)\\
&= [(c-a)F(a-)F_1-(c-a-bx)F F_1]\\
&\hspace{.3in}-[(c-a)F_1(a-)F-(c-a-b(1-x))F_1F]\\
&\hspace{.3in}+(1-2x)FF_1\\
&= (c-a)[F(a-)F_1 -F_1(a-)F] +(1-2x)(1-b)FF_1.
\end{align*}

Since $F_1(a-)F$ is increasing on $(0,1)$, it follows that $f'(x)$ is positive on $(0,1/2)$ and negative on $(1/2,1)$.
\end{proof}
\vspace{.2in}

Particularly interesting hypergeometric functions are the {\it complete elliptic integrals of the first kind}, defined by
$$
\K(x) = \frac{\pi}{2}F\left(\frac 12,\frac 12;1;x^2\right),\ \ \K'(x)=\K(x'),
$$
for $x\in (0,1)$.  Theorem \ref{th4} has the following application to these elliptic integrals.
\vspace{.2in}

\begin{corollary}\label{Sug}
Let $f(x) = x^2x'^2\K(x)\K(x').$  Then  $f(x)$ is increasing on $(0,1/\sqrt 2]$ and decreasing on $[\sqrt{1/2},1)$, with
$$
\max_{0<x<1}f(x) = f\left(\sqrt{1/2}\right) = \frac 14\left(\K\left(\sqrt{1/2}\right)\right)^2 = 0.859398\ldots .
$$
\end{corollary}
\begin{proof}
This follows from Theorem \ref{th4}, if we take $a=b=1/2, c=1$, and replace $x$ by $x^2$.
\end{proof}
\vspace{.2in}

\begin{Pf14}
\end{Pf14}
Here $a_n$, the coefficient of $x^n$ in the hypergeometric series for $F(x)$, is $(a,n)(b,n)/[(c,n)n!]$, so that
$$
n a_n = (a,n)(b,n) /[(c,n) (n-1)!].
$$
This is increasing  if and only if
$n(a+b-c) + ab > 0$, which is true if $a+b \ge c$. Next,
$$
T_n(f(x),1/\sqrt{1-x}) = \frac{n!a_n}{(1/2,n)} = \frac{2^n(a,n)(b,n)}{(c,n)\cdot 1\cdot 3\cdots (2n-1)},
$$
which is decreasing if and only if $2n(a+b-c-1/2) + (2ab -c) < 0$, which is satisfied if
$a+b-1/2 < c$ and $2ab \le c$.  Hence, the assertion follows from Theorem \ref{th2}(7).
$\square$
\vspace{.2in}

\begin{theorem} (cf. \cite[(1.12) and Remark 1.13]{8}  Let $F(x)$ denote the hypergeometric function $F(a,b;a+b;x)$, with $a,b \in (0,1]$ and $|x| < 1$. Then $1/F(x)$ is concave on $(0,\infty )$. In particular,
$$
F\left(\frac{x+y}2\right) \le \frac{2F(x)F(y)}{F(x)+F(y)},
$$
for all $x, y \in (0,1)$, with equality if and only if $x=y$.
\end{theorem}
\begin{proof}
In this case $c = a+b$, and $c-2ab = a(1-b) +b(1-a) \ge 0$, so that the assertion follows from Theorem \ref{newth2}.
\end{proof}
\vspace{.2in}

The next result improves \cite[Theorem 1.25]{8}.
\vspace{.2in}

\begin{theorem}\label{Bessel}
Let $f(x)=\sum_{n=0}^{\infty}b_nx^n$, where $b_n = (-c/4)^n/[n! (k,n)]$, $k = p + (b+1)/2$, as in \cite{8}, be the generalized-normalized Bessel function of the first kind of order $p$. Let $c < 0$, $k > 0$, and $R > 0$. Then
\begin{enumerate}

\item[(1)] $\log f(Re^{-t}))$ is convex on $(0,\infty )$, so that
$f(\sqrt{xy})\le \sqrt{f(x)f(y)}$ for all $x,y\in (0,\infty )$, with equality if and only if  $x = y$.
\item[(2)] $\log f(x)$ is concave on $(0,\infty )$, so that $\sqrt{f(x)f(y)}\le f((x+y)/2)$ for all $x,y\in (0,\infty )$, with equality if and only if $x=y$.
\item[(3)] $f(x)$ is convex on $(0,\infty )$, so that
$$
f((x+y)/2) \le (1/2)(f(x)+f(y))
$$
for all $x,y \in (0,\infty )$, with equality if and only if $x=y$.
\item[(4)] If $k > -1 - cR/4$, then $f(R(1-e^{-t}))$ is concave on $(0,\infty )$, so that $(f(x)+f(y))/2 \le f(R-\sqrt{(R-x)(R-y)})$ for all $x,y\in (0,R)$, with equality if and only if $x=y$.
\end{enumerate}
\end{theorem}
\begin{proof} By the ratio test, the radius of convergence of the series for $f(x)$ is $\infty$.
\begin{enumerate}

\item[(1)] This follows from Theorem \ref{th2}(7) and Lemma \ref{lem1}.
\item[(2)] $T_n(f'(x),f(x)) =(n+1)b_{n+1}/b_n =(-c/4)/(k+n)$, which is decreasing; hence the result follows from Theorem \ref{th2}(2).
\item[(3)] This is obvious, since $b_n > 0$ for all $n$.
\item[(4)] Since
$$
T_n(f'(x),1/(R-x)) = (n+1)b_{n+1}R^{n+1} = (-cR/4)^{n+1}/[n!(k,n+1)]
$$
and $k > -1 -cR/4,$ {\rm we have}
$$
\frac{T_{n+1}}{T_n} = (-cR/4)/[(n+1)(k+n+1)] < 1.
$$
Hence, the result follows from Theorem \ref{th2}(6).
\end{enumerate}
\end{proof}
\vspace{.2in}

\begin{remark}
For $0 < x < y$, let $y/x = \exp (2 \sqrt t)$, $t \in (0,\infty )$, and let
$$
G(x,y) = \sqrt{xy},\ L(x,y) = \frac{y-x}{\log y - \log x},\ {\rm and}\ A(x,y) = \frac{x+y}{2},
$$
denote the {\it Geometric Mean}, {\it Logarithmic Mean}, and {\it Arithmetic Mean} of $x$ and $y$, respectively. Then
\begin{equation*}
\frac{L(x,y)}{G(x,y)} = \frac{\sinh (\sqrt t)}{\sqrt t} = \sum_{n= 0}^{\infty }\frac{t^n}{(2n+1)!},
\end{equation*}
\begin{equation*}
\frac{A(x,y)}{G(x,y)} =  \cosh (\sqrt t) = \sum_{n= 0}^{\infty }\frac{t^n}{(2n)!}.
\end{equation*}
Hence, it will be interesting to study the convexity properties of these two functions.
\end{remark}
\vspace{.2in}

\begin{corollary}\label{corR610} (cf. \cite [Corollary 1.26]{8})
\begin{enumerate}

\item[(1)] \begin{align*}
\cosh (\sqrt{xy}) &\le \sqrt {(\cosh x)(\cosh y)}\\
&\le \cosh \left(\sqrt {(x^2+y^2)/2}\right)\\
&\le \frac 12(\cosh x + \cosh y)
\end{align*}
for all $x,y \in (0,\infty )$, with equality if and only if $x=y$.
\item[(2)]
$$
\frac 12(\cosh x + \cosh y) \le \cosh\left(\sqrt{R-\sqrt{(R-x^2)(R-y^2)}}\right)
$$
for $0<R<6$ and all $x,y\in(0,\sqrt R)$, with equality if and only if $x=y$.
\item[(3)]
\begin{align*}
\frac {\sinh (\sqrt {xy})}{\sqrt {xy}} &\le \sqrt {\frac{\sinh x}{x}\cdot \frac{\sinh y}{y}}\\
&\le \frac{\sinh\sqrt{\frac 12(x^2+y^2)}}{\sqrt{\frac 12(x^2+y^2)}}\\
&\le \frac 12\left(  \frac{\sinh x}{x}  +\frac{\sinh y}{y}\right)
\end{align*}
for all $x,y\in (0,\infty )$, with equality if and only if $x=y$.
\item[(4)]
$$
\frac 12\left(\frac{\sinh x}{x}+\frac{\sinh y}{y}\right)\le \frac{\sinh\left(\sqrt{R-\sqrt{(R-x^2)(R-y^2)}}\right)}{\sqrt{R-\sqrt{(R-x^2)(R-y^2)}}}
$$
for $0<R<10$ and all $x,y\in(0,\sqrt R)$, with equality if and only if $x=y$.
\end{enumerate}
\end{corollary}
\begin{proof}
\begin{enumerate}

\item[(1)] In Theorem \ref{Bessel}, let $b= 1$, $c = - 1$, and $p = -1/2$.
Then $f(x^2) = \cosh x$, and the result follows from Theorem \ref{Bessel}(1),(2),(3) if we replace $x$ and $y$ by $x^2$ and $y^2$, respectively.
\item[(2)]  In Theorem \ref{Bessel}, take $b= 1$, $c = - 1$, and $p = -1/2$.  Then $k=1/2 > -1 + R/4$ if and only if $R<6$.  So the result follows from Theorem \ref{Bessel}(4) if we replace $x$ and $y$ by $x^2$ and $y^2$, respectively.
\item[(3)]  In Theorem \ref{Bessel}, take $b= 1$, $c = -1$, and $p = 1/2$. Then $f(x^2) = (\sinh x)/x$, and the result again follows from Theorem \ref{Bessel}(1),(2),(3), if we replace $x$ and $y$ by $x^2$ and $y^2$, respectively.
\item[(4)]  In Theorem \ref{Bessel}, take $b= 1$, $c = -1$, and $p = 1/2$. Then $k=3/2 > -1 + R/4$ if and only if $R<10$.  So the result follows from Theorem \ref{Bessel}(4) if we replace $x$ and $y$ by $x^2$ and $y^2$, respectively.
\end{enumerate}
\end{proof}
\vspace{.2in}

\begin{remark}
\begin{enumerate}

\item[(1)] The condition $ab/(a+b+1) < c$ in Theorem \ref{th3} cannot be removed. For example, let $f_3(x) = F(3,3;1;x)$ and  $g_3(x) = (d/dx)\log f_3(x).$ According to \cite[p. 484, \#310]{15}, 
$f_3(x) = (1+4x + x^2)/(1-x)^5.$  So 
$$
g_3(x) = (4+2x)/(1+4x+x^2)+5/(1-x).
$$
Clearly $g_3(0) = 9$ and $g_3(0.1)= 8.534\ldots < g_3(0).$
Thus  $f_3$  is not log-convex.  Note that in this example $ab/(a+b+1) = 9/7 > 1 = c$.
\item[(2)] More generally, for $n\ge 3$ let  $f_n(x)=F(n,n;1;x)$ and  $g_n(x)=(d/dx)\log f_n(x)$.  Note that  $ab/(a+b+1) = n^2/(2n+1) > 1$ for $n \ge 3$.  It follows from \cite[Theorem 20, p. 60 and (2), p. 166]{16} that
\addtocounter{equation}{11}
\begin{equation}\label{Fndef}
F(n,n;1;x)=\frac{F(1-n,n;1;\frac x{x-1})}{(1-x)^n}
=\frac{P_{n-1}(y)}{(1-x)^n},
\end{equation}
$-1 < x < 1/2$, where $y=(1+x)/(1-x)$ and  $P_n$ is the Legendre polynomial of degree $n$. These polynomials satisfy the differential recurrence relations \cite[(2), p. 159 and (4), p. 160]{16}
\begin{equation}\label{diffrecur1}
xP_n'(x)=nP_n(x)+P_{n-1}'(x)
\end{equation}
and
\begin{equation}\label{diffrecur2}
P_n''(x)=xP_{n-1}''(x)+(n+1)P_{n-1}'(x).
\end{equation}
To demonstrate that $f_n(x)$ is not log-convex, it will be sufficient to show that $g_n'(0) < 0$, since then $g_n(x)$ cannot be increasing on $(0,1)$.

	Since $P_n(1)=1$ \cite[(5), p. 158]{16}, by (\ref{diffrecur1}) we have $P_n'(1) = n + P_{n-1}'(1).$  By induction and the fact that  $P_1(x)=x$ \cite[p. 160]{16}, we then have
\begin{equation}\label{pn'1}
P_n'(1) = \frac{n(n+1)}{2}.
\end{equation}

	Next, by (\ref{Fndef}),
\begin{equation}\label{gnx}
g_n(x)=\frac{n}{1-x}+\frac{2}{(1-x)^2}\frac{P_{n-1}'(y)}{P_{n-1}(y)}.
\end{equation}
By (\ref{gnx}) and (\ref{pn'1}) we then have  $g_n(0) = n+n(n-1)=n^2,$  in agreement with part (1).  Further,
\begin{align*}
g_n'(x)&=\frac{n}{(1-x)^2}+\frac{4}{(1-x)^3}\frac{P_{n-1}'(y)}{P_{n-1}(y)}\\
&+\frac{4}{(1-x)^4}\frac{P_{n-1}(y)P_{n-1}''(y)-[P_{n-1}'(y)]^2}{[P_{n-1}(y)]^2}.
\end{align*}
With $x=0$, from (\ref{pn'1}) we have
\begin{align*}
g_n'(0) =n+2(n-1)n-(n-1)^2n^2+4P_{n-1}''(1).
\end{align*}
By (\ref{diffrecur2}) with $x=1$ and $P_1''(1)=0$ we have, by induction,
$$
4P_{n-1}''(1)=\frac 12(n-1)^2n^2-(n-1)n.
$$
Thus
$$
g_n'(0) = -\frac 12n^4 + n^3 + \frac 12n^2,
$$
and it is easy to check that this is negative for $n\ge 3$.  In particular, $g_3'(0) = -9$, as a direct computation in part (1) shows.
\item[(3)]  A simple example for which the theorem holds is 
$$
f(x)=F\left(\frac 14,\frac 34;\frac 32;x\right) = \left[\frac{2}{1+\sqrt{1-x}}\right]^{1/2}
$$
(see \cite[p. 472, \#65]{15}), where  $ab/(a+b+1) = 3/32 < 3/2 = c$.  In this case it is easy to see that $f$ is logarithmically convex on $(0,1)$, since we have $(d/dx)\log f(x)= [4(1-x+\sqrt{1-x})]^{-1}$,
which is clearly increasing on $(0,1)$.
\item[(4)] Computer experiments show that in Corollary \ref{corR610}(2) the bound $R<6$ cannot be replaced by $R<7$ and that in Corollary \ref{corR610}(4) the bound $R<10$ cannot be replaced by $R<11\,.$
\end{enumerate}
\end{remark}
\vspace{.2in}
\addtocounter{lemma}{5}

\begin{Pf15}
\end{Pf15}
By Theorem \ref{th2}(1),(5) we have
$$
f(\sqrt{xy}) \le \sqrt{f(x)f(y)} \le f(R-\sqrt{(R-x)(R-y)})
$$
for all $x,y \in (0,R)$, with equality if and only if $x=y$. If we change $x,y$ to (i) $x^2/R$, $y^2/R$ and (ii) $R-x^2/R$, $R-y^2/R$, respectively, then
\begin{equation*}
\begin{split}
f(xy/R) &\le \sqrt{f(x^2/R)f(y^2/R)}\\
&\le f(R-\sqrt{(R-x^2/R)(R-y^2/R)})
\end{split}
\end{equation*}
and
\begin{equation*}
\begin{split}
f(\sqrt{(R-x^2/R)(R-y^2/R)}) &\le \sqrt{f(R-x^2/R) f(R-y^2/R)}\\
&\le f(R-xy/R),
\end{split}
\end{equation*}
for all $x,y \in (0, R)$, with equality  if and only if $x=y$. The result follows if we divide the second chain of inequalities by the first.
$\square$
\vspace{.2in}

We may extend some of the previous results on log-convexity to the generalized hypergeometric function, which is defined as follows.  For non-negative integers $p, q$, let $a_1,\ldots a_p, b_1,\ldots ,b_q$ be positive numbers.  The {\it generalized hypergeometric function} is defined on $(-1,1)$ as
$$
F(x)={}_pF_q(x) = {}_pF_q(a_1,\ldots,a_p;b_1,\ldots b_q;x)=1+\sum_{n=1}^{\infty}\frac{\prod_{k=1}^p(a_k,n)}{\prod_{k=1}^q(b_k,n)}\frac{x^n}{n!},\leqno
$$
where no denominator parameter  $b_k$ is zero or a negative integer.  
\vspace{.2in}

\begin{theorem}
\begin{enumerate}

\item[(1)] If $p=q = 0$, then $F(x) = e^x$, which is trivially log-convex.
\item[(2)] Let $p = q \ge 1$. If $a_k \le b_k$ for each $k$, with at least one strict inequality, then $F$ is strictly log-convex on $(0,1)$.  If $a_k \ge b_k$, with at least one strict inequality, then $F$ is strictly log-concave on $(0,1)$.
\item[(3)] If $p > q$ and $a_k \le b_k$, with at least one strict inequality, for $k = 1,2,\ldots ,q$, then $F$ is strictly log-convex on $(0,1)$.
\item[(4)] If $1\le p < q$ and $a_k \ge b_k$, with at least one strict inequality, for $k = 1,2,\ldots ,p$, then $F$ is strictly log-concave on $(0,1)$.
\item[(5)] If $p = 0$, and $q\ge 1$, then $F$ is log-concave.
\end{enumerate}
\end{theorem}
\begin{proof}
For (1), $F(x) = \sum_{n=0}^{\infty }x^n/n! = e^x$, hence the result.

In case (2),
$$
T_n(F'(x),F(x)) = \frac{B}{A}\frac{a_1 +n}{b_1 +n}\cdots \frac{a_p +n}{b_p +n},
$$
where $A = a_1\cdots a_p$ and
$B = b_1\cdots b_p$. Clearly, a ratio of the form $(a+n)/(b+n)$ is increasing or decreasing in $n$ according as $a < b$ or $a > b$. Hence, $F'(x)/F(x)$ is increasing or decreasing as asserted, and the result follows.

(3) As in case (2), if  $p > q \ge 1$, each ratio of the form $(a_k + n)/(b_k + n)$ is increasing, with at least one strictly, hence so is $F'(x)/F(x)$.  Next, if $q=0$ and $p> 0$, then $T_n(F'(x),F(x)) = (a_1+n)(a_2+n)\cdots (a_p+n)$, which is clearly increasing, so that $F'(x)/F(x)$ is also increasing on $(0,1)$.  Thus $F$ is log-convex.

(4) As in case (2), each ratio of the form $(a_k + n)/(b_k + n)$ is decreasing, with at least one strictly, hence so is $F'(x)/F(x)$.

(5) Here, $T_n(F'(x),F(x)) = 1/[(n+b_1)(n+b_2)\cdots (n+b_q)]$, which is clearly decreasing.
\end{proof}
\vspace{.2in}

3.19. {\it Open problem.} The results of this paper give sufficient
conditions, in terms of the Maclaurin coefficients,
for certain functional inequalities to hold.  What can be said
about the necessary conditions?  Determine functions that
satisfy these inequalities as equalities. 
\vspace{.2in}

{\bf Acknowledgments.}  The authors wish to thank T. Koornwinder for calling our attention to formula (\ref{Fndef}) above in private correspondence.

\vspace{.4in}

\noindent ANDERSON:

\noindent Department of Mathematics

\noindent Michigan State University

\noindent East Lansing, MI 48824, USA

\noindent email: ~~{\tt anderson@math.msu.edu}

\noindent FAX: +1-517-432-1562
\vspace{.2in}

\noindent VAMANAMURTHY:

\noindent Department of Mathematics

\noindent University of Auckland

\noindent Auckland, NEW ZEALAND

\noindent email: ~~{\tt vamanamu@math.auckland.nz}

\noindent FAX: +649-373-7457
\vspace{.2in}

\noindent VUORINEN:

\noindent Department of Mathematics

\noindent FIN-00014, University of Turku

\noindent FINLAND

\noindent e-mail: ~~{\tt vuorinen@utu.fi}

\noindent FAX: +358-2-3336595

\end{document}